%%%This is version of June 24,2007 

\magnification=1200
\input amstex
\input amsppt.sty
\input epsf
\TagsOnRight \openup 3pt
\pagewidth{16.50truecm}
\pageheight{23.00truecm}
\overfullrule=0pt

\def\al{\alpha}

\def\be{\beta}
\def\ep{\varepsilon}

\def\la{\lambda}

\def\ga{\gamma}

\def\cL{\Cal L}

\def\limt{\lim_{t\to\infty}}

\def\limn{\lim_{n\to\infty}}
\def\lims{\lim_{s\to\infty}}

\def\wt{\widetilde}
\def\wh{\widehat}

%Toegevoegd door Rob:

\redefine\l{\ell}

\def\sobre#1#2{\lower 1ex \hbox{ $#1 \atop #2 $ } }
\def\bajo#1#2{\raise 1ex \hbox{ $#1 \atop #2 $ } }

\def\today{\rightline{\ifcase\month\or
  January\or February\or March\or April\or May\or June\or
  July\or August\or September\or October\or November\or December\fi
  \space\number\day,\space\number\year}}
%\vsize 17.7cm  %esta linea es para ver una hoja entera en la pantalla
%\scrollmode
%\baselineskip 15pt
\overfullrule=0pt
%%%%%%%%%%%%%%%%%%
%%%%%%%%%%%%%%%%%%%%%%
%%%%%%%%%%%%
%%% numero del paragrafo e della prima formula dichiarando
%%% \numsec=... \numfor=...  (brevetto Eckmannn).
%%% Si possono citare formule seguenti; le corrispondenze fra nomi
%%% simbolici e numeri effettivi sono memorizzate nel file \jobname.aux, che
%%% viene letto all'inizio, se gia' presente. E' possibile citare anche
%%% formule che appaiono in altri file, purche' sia presente il
%%% corrispondente file .aux; basta includere all'inizio l'istruzione
%%%           \include{nomefile}
%%%
%%%%%%%%%%%%%%%%%%%%%%%%%%%%%%%%%%%%%%%%%%%%%%%%%%%%%%%%%%%%%%%

\global\newcount\numsec
\global\newcount\numfor
\global\newcount\numtheo
\global\advance\numtheo by 1

\def\senondefinito#1{\expandafter\ifx\csname#1\endcsname\relax}

\def\SIA #1,#2,#3 {\senondefinito{#1#2}%
\expandafter\xdef\csname #1#2\endcsname{#3}\else
\write16{???? ma #1,#2 e' gia' stato definito !!!!} \fi}

\def\etichetta(#1){(\veroparagrafo.\veraformula)%
\SIA e,#1,(\veroparagrafo.\veraformula) %
\global\advance\numfor by 1%
\write15{\string\FU (#1){\equ(#1)}}%
\write16{ EQ #1 ==> \equ(#1) }}

\def\letichetta(#1){\veroparagrafo.\verotheo
\SIA e,#1,{\veroparagrafo.\verotheo}
\global\advance\numtheo by 1
 \write15{\string\FU (#1){\equ(#1)}}
 \write16{ Sta \equ(#1) == #1 }}

\def\tetichetta(#1){\veroparagrafo.\veraformula %%%%copy four lines
\SIA e,#1,{(\veroparagrafo.\veraformula)}
\global\advance\numfor by 1
\write15{\string\FU (#1){\equ(#1)}}
\write16{ tag #1 ==> \equ(#1)}}

\def\FU(#1)#2{\SIA fu,#1,#2 }

\def\etichettaa(#1){(A\veroparagrafo.\veraformula)%
\SIA e,#1,(A\veroparagrafo.\veraformula) %
\global\advance\numfor by 1%
\write15{\string\FU (#1){\equ(#1)}}%
\write16{ EQ #1 ==> \equ(#1) }}

\def\BOZZA{
\def\alato(##1){%
 {\rlap{\kern-\hsize\kern-1.4truecm{$\scriptstyle##1$}}}}%
\def\aolado(##1){%
 {%\vtop to \profonditastruttura
{%\baselineskip
 %\profonditastruttura\vss
 \rlap{\kern-1.4truecm{$\scriptstyle##1$}}}}}
}

\def\alato(#1){}
\def\aolado(#1){}

\def\veroparagrafo{\number\numsec}
\def\veraformula{\number\numfor}
\def\verotheo{\number\numtheo}

\def\Eq(#1){\eqno{\etichetta(#1)\alato(#1)}}
\def\eq(#1){\etichetta(#1)\alato(#1)}
\def\teq(#1){\tag{\tetichetta(#1)\hskip-1.6truemm\alato(#1)}}
%%%%%for \tag
\def\Eqa(#1){\eqno{\etichettaa(#1)\alato(#1)}}
\def\eqa(#1){\etichettaa(#1)\alato(#1)}
\def\eqv(#1){\senondefinito{fu#1}$\clubsuit$#1
\write16{#1 non e' (ancora) definito}%
\else\csname fu#1\endcsname\fi}
\def\equ(#1){\senondefinito{e#1}\eqv(#1)\else\csname e#1\endcsname\fi}

%next six lines by paf (no responsibilities taken)
\def\Lemma(#1){\aolado(#1)Lemma \letichetta(#1)}%
\def\Theorem(#1){{\aolado(#1)Theorem \letichetta(#1)}}%
\def\Proposition(#1){\aolado(#1){Proposition \letichetta(#1)}}%
\def\Corollary(#1){{\aolado(#1)Corollary \letichetta(#1)}}%
\def\Remark(#1){{\noindent\aolado(#1){\bf Remark \letichetta(#1)}}}%
\def\Definition(#1){{\noindent\aolado(#1){\bf Definition \letichetta(#1)}}}
%$\!\!$\hskip-1.6truemm
\def\Example(#1){\aolado(#1) Example \letichetta(#1)$\!\!$\hskip-1.6truemm}

\def\Assumptions(#1){{\noindent\aolado(#1){\bf Assumptions \letichetta(#1)}}}
\def\include#1{
\openin13=#1.aux \ifeof13 \relax \else
\input #1.aux \closein13 \fi}

\openin14=\jobname.aux \ifeof14 \relax \else
\input \jobname.aux \closein14 \fi
\openout15=\jobname.aux

%%% Local Variables:
%%% mode: plain-tex
%%% TeX-master: t
%%% End:

%%%%%%%%References macros start here%%%%%%%%%%
%
\catcode`\X=12\catcode`\@=11
\def\n@wcount{\alloc@0\count\countdef\insc@unt}
\def\n@wwrite{\alloc@7\write\chardef\sixt@@n}
\def\n@wread{\alloc@6\read\chardef\sixt@@n}
\def\crossrefs#1{\ifx\alltgs#1\let\tr@ce=\alltgs\else\def\tr@ce{#1,}\fi
   \n@wwrite\cit@tionsout\openout\cit@tionsout=\jobname.cit
   \write\cit@tionsout{\tr@ce}\expandafter\setfl@gs\tr@ce,}
\def\setfl@gs#1,{\def\@{#1}\ifx\@\empty\let\next=\relax
   \else\let\next=\setfl@gs\expandafter\xdef
   \csname#1tr@cetrue\endcsname{}\fi\next}
\newcount\sectno\sectno=0\newcount\subsectno\subsectno=0\def\r@s@t{\relax}
\def\resetall{\global\advance\sectno by 1\subsectno=0
  \gdef\firstpart{\number\sectno}\r@s@t}
\def\resetsub{\global\advance\subsectno by 1
   \gdef\firstpart{\number\sectno.\number\subsectno}\r@s@t}
\def\v@idline{\par}\def\firstpart{\number\sectno}
\def\l@c@l#1X{\firstpart.#1}\def\gl@b@l#1X{#1}\def\t@d@l#1X{{}}
\def\m@ketag#1#2{\expandafter\n@wcount\csname#2tagno\endcsname
     \csname#2tagno\endcsname=0\let\tail=\alltgs\xdef\alltgs{\tail#2,}%
  \ifx#1\l@c@l\let\tail=\r@s@t\xdef\r@s@t{\csname#2tagno\endcsname=0\tail}\fi
   \expandafter\gdef\csname#2cite\endcsname##1{\expandafter
 %the following line was replaced by the subseqent one, DNA 7/6/89
  %  \ifx\csname#2tag##1\endcsname\relax?\else\csname#2tag##1\endcsname\fi
     \ifx\csname#2tag##1\endcsname\relax?\else{\rm\csname#2tag##1\endcsname}\fi
    \expandafter\ifx\csname#2tr@cetrue\endcsname\relax\else
     \write\cit@tionsout{#2tag ##1 cited on page \folio.}\fi}%
   \expandafter\gdef\csname#2page\endcsname##1{\expandafter
     \ifx\csname#2page##1\endcsname\relax?\else\csname#2page##1\endcsname\fi
     \expandafter\ifx\csname#2tr@cetrue\endcsname\relax\else
     \write\cit@tionsout{#2tag ##1 cited on page \folio.}\fi}%
   \expandafter\gdef\csname#2tag\endcsname##1{\global\advance
     \csname#2tagno\endcsname by 1%
   \expandafter\ifx\csname#2check##1\endcsname\relax\else%
\fi%      \immediate\write16{Warning: #2tag ##1 used more than once.}\fi
   \expandafter\xdef\csname#2check##1\endcsname{}%
   \expandafter\xdef\csname#2tag##1\endcsname
     {#1\number\csname#2tagno\endcsnameX}%
   \write\t@gsout{#2tag ##1 assigned number \csname#2tag##1\endcsname\space
      on page \number\count0.}%
   \csname#2tag##1\endcsname}}%
\def\m@kecs #1tag #2 assigned number #3 on page #4.%
   {\expandafter\gdef\csname#1tag#2\endcsname{#3}
   \expandafter\gdef\csname#1page#2\endcsname{#4}}
\def\re@der{\ifeof\t@gsin\let\next=\relax\else
    \read\t@gsin to\t@gline\ifx\t@gline\v@idline\else
    \expandafter\m@kecs \t@gline\fi\let \next=\re@der\fi\next}
\def\t@gs#1{\def\alltgs{}\m@ketag#1e\m@ketag#1s\m@ketag\t@d@l p
    \m@ketag\gl@b@l r \n@wread\t@gsin\openin\t@gsin=\jobname.tgs \re@der
    \closein\t@gsin\n@wwrite\t@gsout\openout\t@gsout=\jobname.tgs }
\outer\def\localtags{\t@gs\l@c@l}
\outer\def\globaltags{\t@gs\gl@b@l}
\outer\def\newlocaltag#1{\m@ketag\l@c@l{#1}}
\outer\def\newglobaltag#1{\m@ketag\gl@b@l{#1}}

\def\t@gsoff#1,{\def\@{#1}\ifx\@\empty\let\next=\relax\else\let\next=\t@gsoff
   \expandafter\gdef\csname#1cite\endcsname{\relax}
   \expandafter\gdef\csname#1page\endcsname##1{?}
   \expandafter\gdef\csname#1tag\endcsname{\relax}\fi\next}
\def\verbatimtags{\let\ift@gs=\iffalse\ifx\alltgs\relax\else
   \expandafter\t@gsoff\alltgs,\fi}
\catcode`\X=11 \catcode`\@=\active \localtags
%%%%%%%%%%%%references macro end here%%%%%%%%%%%%%%%

%%%%%%%%%%%%%%%%%%%%%%%%%%%%%%%%%
%%%%%%%%%%%%%%%%%%%%%%%%%%%%%%%%%%
%%%%%%%%%%%%

%%%%%%%%%%%%%%%%%% Numerazione verso il futuro ed eventuali paragrafi
%%%%%%%      precedenti non inseriti nel file da compilare
\def\include#1{
\openin13=#1.aux \ifeof13 \relax \else
\input #1.aux \closein13 \fi}
\openin14=\jobname.aux \ifeof14 \relax \else
\input \jobname.aux \closein14 \fi
%\openout15=\jobname.aux
%%%%%%%%%%%%%%%%%%%%%%%%%%%%

\BOZZA
\topmatter
\title 
A problem in last-passage percolation
\endtitle
\author 
Harry Kesten and Vladas Sidoravicius
\endauthor

\abstract
Let $\{X(v), v \in \Bbb Z^d \times \Bbb Z_+\}$ be an i.i.d. family of
random variables such that $P\{X(v)= e^b\}=1-P\{X(v)= 1\} = p$ for some
$b>0$. We consider paths $\pi \subset \Bbb Z^d \times \Bbb Z_+$ starting 
at the origin and with the last coordinate increasing along the path, and 
of length $n$. Define for such paths $W(\pi) = \text{number of vertices $\pi_i, 1 \le i \le n$, with }X(\pi_i) =
e^b$. Finally let  $N_n(\al) = \text{number of paths $\pi$ of length $n$ starting at $\pi_0
= \bold 0$ and with $W(\pi) \ge \al n$.}$ We establish several properties
of $\lim_{n \to \infty} [N_n]^{1/n}$.
\endabstract

\endtopmatter
\document

\numsec1
\numfor1
\subhead
1. Statement of the problem
\endsubhead
The study of the free energy of a directed polymer in random environment
suggested the problems of this paper to us. Here we consider a site 
version of semi-oriented first-passage percolation. To
be more precise we take for $\cL$ 
the graph $\Bbb Z^d\times \Bbb Z_+$ with the last coordinate 
oriented in the standard way. 
A vertex $v \in \Bbb Z_+^d$ has an edge
to $v\pm e_i +e_{d+1}$ for $1 \le i\le d$,
and there are no other outgoing edges from $v$. 
Here and in the sequel $e_i$ stands for the $i$-th coordinate vector.
We shall use the symbol $\bold 0$ for the origin in $\Bbb Z^d$, as
well as for the corresponding vertex of $\cL$.
For $v=(v_1, \dots,v_d)$ a vertex of $\cL$ or of $\Bbb Z^d$, $\|v\|$
will be the $\l_1$-norm of $v$, i.e., $\|v\| = \sum_{i=1}^d |v_i|$.
We will call a path on $\cL$ {\it semi-oriented} 
and we will say that we are dealing
with the semi-oriented case.

Our arguments can also be carried out in a the fully oriented  case in which 
 $\cL$ is replaced by the graph $\Bbb Z^{d+1}_+$
with an edge from $v$ to $v+e_i$ for $v \in \Bbb
Z^{d+1}_+$ and $1 \le i \le d+1$. However, we shall not mention the
latter case anymore in these notes.

We assign to each $v \in \cL$ a random weight $X(v)$. The $\{X(v):v
\in \cL\}$ are taken i.i.d. with the common distribution 
$$
P\{X(v) = e^b\} = p, \quad P\{X(v) = 1\} =1 -p,
$$
for some $b > 0, 0< p < 1$. 
Nothing interesting happens when $p=0$ or 1, so we exclude these 
values for $p$. 
For an oriented path $\pi =(\pi_0, \pi_1,
\dots, \pi_s)$ on $\cL$ of length $s$ we define
$$
W(\pi) = \text{number of vertices $\pi_i, 1 \le i \le s$, with }X(\pi_i) =
e^b.
$$
(Note that $X(\pi_0)$ does not contribute to $W(\pi)$.)  
We further define for $0 \le \al \le 1-p$ 
$$
N_s(\al) = \text{number of paths $\pi$ of length $s$ starting at $\pi_0
= \bold 0$ and with $W(\pi) \ge \al s$.}
$$
We are interested in these notes in the behavior of $N_s(\al)$ for large
$s$ and different $\al$.

We have been informed that related problems have been studied by 
\cite{CPV}.

The first lemma is an exponential bound for $P\{N_s(\al) = 0\}$ for
certain $\al$, as $s
\to \infty$. Basically this comes from \cite{GK}, but the
oriented case considered here is simpler than the unoriented case of
\cite{GK}. See also \cite{CMS}
\proclaim{Lemma 1} The limit
$$
M = M(p) := \lims \max_{\pi_0 = \bold 0, |\pi|=s} \frac 1{|\pi|}W(\pi)
\teq(1.1)
$$
exists and is constant a.s. If $p >0$, then also $M > 0$.
(Here $|\pi|=s$ in the max means that we take the maximum over all
oriented paths of length $s$ which start at $\bold 0$.) 
Moreover, for any $\ep >0$ there exist constants $0 < C_i < \infty$
for which
$$
P\{N_t(M(p)-\ep) = 0\} = P\Big\{\max_{\pi_0 = \bold 0, |\pi|=t} 
\frac {W(\pi)}t < M(p) -\ep\Big\} \le C_1e^{-C_2t}, t \ge 0.
\teq(1.2)
$$
\endproclaim
\demo{Proof}
In the sequel a path will always mean an oriented path on $\cL$.
However, a path does not have to start at $\pi_0$ at time 0. We will call 
the sequence $(\pi_j, \dots, \pi_{j+t})$ a {\it path starting at $\pi_j$ at
time $s$} and of length $t$ if $\|\pi_j\| = s$ and
there is an oriented edge of $\cL$
from $\pi_i$ to $\pi_{i+1}$ for $j \le i < j+t$.

The limit $M$ exists and is a.s. constant by \cite{GK}.
In the oriented case considered here this was proven in an easier way
in \cite{CMS} by an
application of Liggett's subadditive ergodic theorem
(\cite{Li}). We merely outline the proof of \cite{CMS}. Define
$$
M_s(x, y) = \max \Sb \pi_0 = x, |\pi| = s\\
\pi_s = y \endSb W(\pi),
$$
$$
M_s(x,*) = \max_y \;M_s(x,y) = \max_{\pi_0 = x, |\pi| = s} W(\pi).
$$
Define further
$$
y(s) = \text{ first vertex $y$ in lexicographical order for which }
M_s(\bold 0,y)= M_s(\bold 0,*).
$$
Then, for $s,t \ge 1$
$$ 
M_{s+t}(\bold 0,*) \ge M_s(\bold 0,*) +M_t(y(s),*) = M_s(\bold 0,
y(s))+ M_t(y(s),*).
\teq(1.3)
$$
Indeed, the left hand side is a maximum over all paths starting at
$\bold 0$ and of length $s+t$, while the right hand side is just a
maximum over paths which start at $\bold 0$ but pass through $y(s)$ at
time $s$ and have length $s+t$. If one sets $M_0(x,y)=0$ for all $x,y$,
then \equ(1.3) remains valid even if $s=0$ or $t=0$.

We note further that if all $X(v)$
with  $\|v\| \le s$ are given, then $y(s)$ is also fixed and $M_t(y(s),*)$
is defined in the same way as $M_t(\bold 0,*)$, but with $X(v)$ replaced by 
$X(v+y(s))$. It follows from this that the conditional distribution of
$M_t(y(s), *)$ given all $X(v)$ with $\|v\| \le s$ is just the  same
as the unconditional distribution of $M_t(\bold 0,*)$, and hence does
not depend on the $X(v)$ with $\|v\| \le s$. Thus, $M_t(y(s),*)$ is
independent of those $X(v)$ and has the distribution  of $M_t(\bold 0,*)$.
These observations allow us to apply Liggett's theorem
(\cite{Li},Theorem VI.2.6) to the variables $X_{s,t}:= M_{t-s}(y(s),
u(s,t))$, where
$$
u(s+t) =  \text{ first vertex $u$ in lexicographical order for which }
M_t(y(s), u) = M_t(y(s),*).
$$
This shows that $M(p)$ exists and is almost surely constant. 
The fact that $M > 0$ is immediate from
$$
M \ge \limt \frac {W(\pi^{(t)})}t,
$$
where $\pi^{(t)}$ is the path which moves along the first
coordinate axis from $\bold 0$ to $(t, 0, \dots,0)$ in $t$ steps.
Indeed 
$$
\frac {W(\pi^{(t)})}t = \frac 1t\sum_{i=1}^t I[X(i,0, \dots,0) = e^b]
$$
and this tends to $p$ by the strong law of large numbers.

Now, to start on the proof of \equ(1.2) note first that the
equality of the first and second member in \equ(1.2) is immediate from
the definitions. Indeed, $N_t(\al) = 0$ means that for all path $\pi$ of
length $t$ and starting at the origin $W(\pi) \le \al t$.
We therefore concentrate on the inequality in \equ(1.2).
Observe that by definition of $W$
$$
\frac 1{|\pi|}W(\pi)\le 1
\teq(1.25)
$$
so that also 
$$
\frac {M_s}s \text{ is bounded and }\lims \frac {E M_s}s = M.
\teq(1.4)
$$
Let $\ep > 0$ be given. One can then fix $s$  such
that $EM_s/s- \ep/2 \ge M(p)-\ep$.
Now define recursively $y_0 = \bold 0, y_1 = y(s)$, 
$$
y_{k+1} = 
\text{ first vertex $y$ in lexicographical order for which
}M_s(y_k,y) = M_s(y_k, *). 
$$
Analogously to \equ(1.3) we then have 
$$
M_{ks+t}(x,*) \ge M_{ks}(x,z)+ M_t(z,*).
$$
This hold for any $z$ and in particlar for any $z$ for which
$M_{ks}(x,z) = \max_v M_{ks}(x,v)$. By iteration,
$$
\align
M_{ks}(\bold 0,*) &\ge M_s(\bold 0, y_1)+M_{(k-1)s}(y_1,*) \ge M_s(\bold
0,y_1)+ M_s(y_1, y_2)+ M_{(k-2)s}(y_2,*)\\
&\ge \dots  
\ge \sum_{j=0}^{k-1} M_s(y_j,y_{j+1}).
\endalign
$$
By the argument given a few lines after \equ(1.3), the random
variables $M_s(y_k,y_{k+1})$ are i.i.d. Moreover, the variables 
$M_s(y_j,y_{j+1}), j \ge 0,$ are bounded (see \equ(1.25). By
exponential bounds for the sum of i.i.d. variables or Bernstein's
inequality (see \cite{CT}, exercise 4.3.14) we have
$$
\align
&P\{M_{ks}(\bold 0,*) \le ks[M -\ep]\}\\
&\le P\Big\{\sum_{j=0}^{k-1} M_s(y_j,y_{j+1}) 
\le k[EM_s(\bold 0,y_1) - \ep s/2]\Big\}
\le C_1 e^{-C_2k}.
\teq(1.8)
\endalign
$$ 
This proves \equ(1.2) for $t$ a multiple of $s$. The extension to
arbitrary positive integers $t$ is an easy monotonicity argument. If
$ks \le t < (k+1)s$ and $\pi$ is a path of length $t$, let $\pi'$ be
the initial piece of length $ks$ of $\pi$. Then 
$N_t(M-2\ep) =0$ implies $W(\pi') \le W(\pi)\le t(M-2\ep) \le
ks(M-\ep)$ for large $k$ and this happens only on a set of 
probability at most $C_1\exp[-C_2ks]$.
\hfill $\blacksquare$
\enddemo

\medskip
The next lemma will help us to formulate a concrete problem. 
\proclaim{Lemma 2} For $0 \le \al  \ne M$
$$
\la(\al) = \la(\al,p) := \limt [N_t(\al)]^{1/t} \text{ exists and is 
constant a.s.}
\teq(1)
$$
\endproclaim
\demo{Proof} This proof uses standard arguments for superconvolutive
sequences. However the assumptions here seem to differ from the usual
ones and we see no way to appeal to a standard theorem such as
\cite{H} for the lemma. We therefore go into some detail.
We break the proof into 3 steps.
\newline
{\bf Step 1.} To begin with, if $\al > M$, then by the fact that the limit in
\equ(1.1) exists we have $\max_{\pi_0= \bold 0,|\pi|=s} W(\pi)/s <\al$
eventually. But this says that $N_s(\al) = 0$ for all large $n$. Thus
\equ(1) with $\la(\al) = 0$ is obvious when $\al > M$.

Next fix an $\al$ with $\al < M$. We shall suppress $\al$ in our
notation for the rest of this proof. In the rest of this step we
define $N_t$ and related quantities and show that they are almost 
superconvolutive. Define
$$
\align
N_t(x) &= N_t(x;\al)\\
& = \text{ number of paths $\pi$ of length $t$ which
start at $x$ and have $W(\pi) \ge \al t$},
\endalign
$$
$$
\align
N_t(x,y) = N_t(x,y;\al) =&\text{ number of paths $\pi$ of length $t$ which
start at $x$}\\
&\phantom {a}\text{and end at $y$ and have $W(\pi) \ge \al t$},
\endalign
$$
and 
$$
N_t(x, *) = \max_y N_t(x,y).
$$
Note that 
$$
N_t(x) = \sum_y N_t(x,y) \text{ and }N_t = N_t(\bold 0).
$$
Accordingly we set
$$
N_t(*) = N_t(\bold 0,*).
$$
Note also that $N_t(x,y)$ can be nonzero only if $y = x+v$ for some $v \in
\cL$ with $\|v\|= t$. There are at most $(t+1)^d$ possible values for
$v$. Thus the max here is really a maximum over at
most $(t+1)^d$ values of $y$. Consequently,
$$
(t+1)^{-d}N_t (x) \le N_t(x,*)  \le N_t(x).
\teq(2)
$$
It follows from this that it suffices for \equ(1) to prove that
$$
 \limt[N_t(*)]^{1/t} \text{ exists and is constant a.s.}
\teq(3)
$$

The advantage of $N_t(*)$ is that it is almost superconvolutive. To make
this precise, we order the vertices of $\cL$ 
lexicograhically. If $N_t > 0$, then also $N_t(*) >0$. In
this case we define
$$
z(t) = \text{ first site $z$ in the lexicographical ordering for which
}N_t(\bold 0, z) = N_t(*).
$$
If $N_t = 0$, then also $N_t(*) = 0$. In this case we take for
$z(t)$ any fixed vertex $z$ of $\cL$ with
$\|z\| = s$. For the sake of definiteness we shall take $z(t)
=(t,0, \dots,0)$.
With these definitions we have  for $s,t \ge 1$
$$
N_{s+t}(*) \ge N_s(*) \cdot N_t(z(s), *).
\teq(6)
$$
This is trivial if $N_s  = 0$, for then also $N_s(*) =
0$. If $N_s > 0$, and hence also $N_s(*) > 0$, then \equ(6)
follows from the fact that $N_{s+t}(*)$
is no smaller than (number of paths
$\pi =(\pi_0, \dots,\pi_s)$ of length $s$ with $\pi_0 = \bold 0,
\pi_s = z(s)$ and $W(\pi) \ge \al s$) times
(number of paths $\wt \pi$ which start at time $s$ at $\pi_s=z(s)$ and are at
time $s+t$ at any fixed vertex $z$, and have $W(\wt \pi) \ge \al t$). 
The maximum over all $z$ of the second factor is just $N_t(z(s),*)$. 
\newline
{\bf Step 2.} In this step we show that
$$
\limt  \frac 1t E \{\log N_t(*)\} \text{ exists and lies in }[0,\log d].
\teq(6.1)
$$
We set
$$
Y_s= Y_s(\al)= [\log N_s(*)]^+,\;Y_{s,t} =[\log N_t(z(s),*)]^+
$$
and
$$
Z_s= Z_s(\al) = s \log (2d) - Y_s,\;
Z_{s,t} = t\log (2d) - Y_{s,t}.
$$
Note that $Y_s$ is at most equal to the logarithm of the number of
paths $\pi$ of length $s$ with $\pi_0 = \bold 0$, i.e., $Y_s \le s\log
(2d)$. Consequently, 
$$
0 \le Z_s \le s \log (2d).
\teq(7.1)
$$
Similarly,
$$
0 \le Z_{s,t} \le t\log (2d).
\teq(7.2)
$$
On the event $A(s,t) := \{N_s > 0, N_t(z(s)) > 0\}$ it holds $N_s(*) \ge 1$ and
$N_t(z(s),*) \ge 1$, so that $Z_s = s \log (2d) - \log N_s(*)$ and $Z_{s,t} =
t\log (2d) - \log N_t(z(s),*)$. The relation \equ(6) therefore shows
that on the event $A(s,t)$ we have
$$
Z_{s+t} \le Z_s + Z_{s,t}.
\teq(7)
$$ 
Off the event $A_{s,t}$ we need to introduce a correction term. We
define
$$
\Psi(s,t) = \Psi(s,t, \al):=I[N_s = 0]Y_{s,t} +I[N_t(z(s)) = 0]Y_s.
\teq(7.7)
$$
It is now easy to see that always
$$
Z_{s+t} \le Z_s + Z_{s,t} + \Psi(s,t);
\teq(8)
$$
in fact, if $N_s = 0$, then $Y_s = 0$ and the right hand side 
equals $(s+t)\log (2d)$. Similarly if $N_t(z(s))= 0$.

We claim that $N_t(z(s))$ is independent of all $X(v)$ with
$\|v\| \le s$ and has the same distribution as $N_t$. In fact,
if we fix all $X(v)$ with $\|v\| \le s$, then also $z(s)$ is
determined, and $N_t(z(s))$ is defined in the same way as $N_t(\bold 0)
= N_t$, but with $X(v)$ replaced by $X(z(s)+v)$. This shows that the
conditional distribution of $N_t(z(s))$ , given all $X(v)$ with $\|v\|
\le s$ is the same as the unconditional distribution of $N_t$, which
proves our claim.
Taking expectations in \equ(8) therefore gives
$$
\align
EZ_{s+t} &\le E Z_s + EZ_{s,t} + E\Psi(s,t)\\ 
&\le  EZ_s +EZ_t +P\{N_s
=0\}t \log (2d)+ P\{N_t = 0\}s \log (2d).
\teq(9.1)
\endalign
$$
Note that all these expectations are finite by virtue of \equ(7.1) and
\equ(7.2). In particular, if $K$ is any  positve integer, and
$s=t=K2^j$, then
$$
\frac 1{K2^{j+1}}EZ_{K2^{j+1}} \le \frac 2{K2^{j+1}} EZ_{K2^j} 
+ \frac {K2^{j+1}\log (2d)}{K2^{j+1}} P\{N_{K2^j} = 0\}.
$$
However, if we take $\ep = M-\al$, then we see from Lemma 1 that
$$
P\{N_t = 0\} \le C_1 \exp[-C_2t],
\teq(10)
$$
whence
$$
\frac 1{K2^{j+1}}EZ_{K2^{j+1}} \le \frac 1{K2^j} EZ_{K2^j} + 
\log (2d) C_1 \exp[-C_2K2^j].
$$
This easily implies 
$$
\limsup_{j \to \infty}  \frac 1{K2^j} EZ_{K2^j} \le \liminf_{j \to 
\infty}  \frac 1{K2^j} EZ_{K2^j},
$$
so that
$$
\ga(K) := \lim_{j \to \infty}  \frac 1{K2^j} EZ_{K2^j} \text{ exists and
lies in }[0,\log (2d)]
\teq(9)
$$
(see \equ(7.1) for the bounds on $\ga$).

Next we will prove that $\ga(K)$ is independent of $K$. 
Let $K, L \ge 1$ be integers and let the dyadic expansion of
$L/K$ be
$$
\frac LK = \sum_{j= - \infty}^n 2^{k_j},
\teq(9.7)
$$ 
where $k_j$ is increasing in $j$, sign$(k_j)$ = sign$(j)$ and $n$ some
finite non-negative integer. The sum over negative $j$ may actually
be finite, but in order to avoid further notation we write sum over
the negative $j$ as starting at $-\infty$.

The expansion \equ(9.7) can also be written as 
$$
L2^\l = K\sum_{j= - \infty}^n 2^{k_j+\l}
$$
for any integer $\l \ge 0$. We shall let $\l \to \infty$ later on, but
for the moment leave it unspecified. Since we shall use $Z_m$ for
somewhat messy  $m$'s we shall write $Z(m)$ instead of $Z_m$ in the
calculations below. Start with an application of \equ(9.1)
with $s+t = L2^\l, s = K2^{k_n+\l}$. Thus we take 
$$
t = L2^\l-K2^{k_n+\l} = K\sum_{j = - \infty}^{n-1} 2^{k_j+\l}.
\teq(9.2)
$$
Since the right hand side is positive,  $t$ is a
positive integer. Taking into account that 
$$
t = \text{ right hand side of \equ(9.2) }\le K2^{k_n + \l} = s,
$$
we obtain
$$
\align
EZ(L2^\l) &\le EZ(K2^{k_n+\l}) +EZ(t) + P\{N_s = 0\}t \log (2d) + P\{N_t
=0\} s\log (2d) \\
&\le EZ(K2^{k_n+\l}) +EZ(t) + C_1t\log (2d) \exp[-C_2 s]\\
&\phantom{MMMMMMMMMMM}+ 
C_1s\log (2d) \exp[-C_2t] \text{ (by \equ(10))}\\
&\le  EZ(K2^{k_n+\l}) +EZ(t) + L2^\l C_1 \exp[-C_2t]\log (2d).
\teq(10.4)
\endalign
$$
Divide both sides of the inequality by $L2^\l$ and let $\l \to \infty$,
and note that $t \to \infty$ as $\l \to \infty$ (see
\equ(9.2)). \equ(9) then shows that
$$
\ga(L) \le \ga(K) \frac KL2^{k_n} + \limsup_{\l \to \infty} \frac
{EZ(t)}{L2^\l}.
$$
We repeat this argument in the following way. Set 
$$
t_r = K\sum_{j = -\infty}^{n-r} 2^{k_j+\l},  r \ge 0,
$$
and apply \equ(9.1) and \equ(10) with $t_r$ for $s+t$ and
$K2^{k_{n-r}+\l}$ for $s$, and consequently $t_{r+1}$ for $t$.
Taking into account that $t_{r+1} \le K2^{k_{n-r}}$ we obtain
$$
EZ(t_r) \le EZ(K2^{k_{n-r}+\l}) + EZ(t_{r+1}) + C_1t_r\log (2d)
\exp[-C_2t_{r+1}],\; r\ge 0.
\teq(10.5)
$$
For $r=0$ this is just \equ(10.4) with $L2^\l$ for $t_0$. This 
time we successively use 
\equ(10.5) for $r=0, 1, \dots, R-1$ before we divide by $L2^\l$, where
$R$ is determined as follows: (i) if the expansion in \equ(9.7) has
only finitely many terms, then we take $R$ such that $2^{k_{n-R}}$ is
the smallest power of 2 appearing in the right hand side of \equ(9.7)
(so that $t_{R+1} = 0$); (ii) if the expansion in \equ(9.7) has
infintely many terms, then we
fix a small number $\eta > 0$ and let
$R = R(\eta)$ be the smallest non-negative integer such that
$$
K\sum_{j=- \infty}^{n-R} 2^{k_j} \le \eta.
\teq(10.6)
$$
Note that $R$ does not depend on $\l$.
We get
$$
\align
EZ(L2^\l) &= EZ(t_0) \le EZ(K2^{k_n+\l})+ EZ(t_1) 
+ C_1t_0 \log (2d) \exp[-C_2t_1]\\
&\le EZ(K2^{k_n+\l})+ EZ(K2^{k_{n-1}+\l}) + EZ(t_2) 
+ C_1t_0 \log (2d) \exp[-C_2t_1]\\
&\phantom{MMMMMMMM}+C_1t_1\log (2d)\exp[-C_2t_2]\le \dots\\
&\le \sum_{r=0}^{R-1} EZ(K2^{k_{n-r}+\l}) + EZ(t_R) + 
C_1 \log (2d)\sum_{r=0}^{R-1} t_r \exp[-C_2t_{r+1}].
\teq(9.8)
\endalign
$$
Now we divide by $L2^\l$ and let $\l \to \infty$. 
Consider first case (i) when the expansion in \equ(9.7) is finite. 
Now recall 
$$
\frac{t_r}{L2^\l}= \frac KL\sum_{j = -\infty}^{n-r}2^{k_j}\le 1.
$$
On the other hand, $t_{r+1} \to \infty$ as $\l \to \infty$ for each
$r \le R-1$.
%% and $[L2^\l]^{-1}t_r$ is bounded for $r \le R$ in this case.
%%Moreover, $t_R = K2^{k_{n-R}+\l}$, still in case (i).
The inequality \equ(9.8)  therefore implies
$$
\frac 1{L2^\l}C_1 \log (2d)\sum_{r=0}^{R-1} t_r \exp[-C_2t_{r+1}] \to 0
$$
and, by virtue of \equ(9) and $t_R = K2^{k_{n-R}+\l}$,
$$
\ga(L) = \lim_{\l \to \infty} \frac {EZ(L2^\l)}{L2^\l} 
\le  \sum_{j=0}^R \ga(K)\frac KL 2^{k_{n-j}} =\ga(K).
\teq(9.9)
$$

Next, in case (ii) we obtain similarly
$$
\ga(L) \le  \sum_{j=0}^{R-1} \ga(K)\frac KL 2^{k_{n-j}}+ \limsup_{\l \to
\infty}\frac {EZ(t_R)}{L2^\l} + \limsup_{\l \to \infty} C_1 \log (2d)
\exp[-C_2t_R].
$$
This time we use that 
$$
\frac 1{L2^\l}EZ(t_R) \le \frac 1{L2^\l}t_R \log (2d) \text{ (by
\equ(7.1)}) \le \frac {\eta \log (2d)}L \text{ (by \equ(10.6))}.
$$
Finally, 
$$
t_R \ge K2^{k_{n-R-1}}2^\l \to \infty \text{ as } \l \to \infty,
$$
because the term $2^{k_{n-R-1}}$ is actually present in \equ(9.7) in case (ii).
Thus in case (ii)
$$
\ga(L) \le  \sum_{j=0}^{R-1} \ga(K)\frac KL 2^{k_{n-j}} + \frac{\eta
\log (2d)}L \le \ga(K)  + \frac{\eta \log (2d)}L.
$$
Since this holds for any $\eta > 0$ we obtain in both cases that $\ga(L)
\le \ga(K)$. By interchanging the roles of $K$ and $L$ we finally 
prove that $\ga(K)$ does not depend on $K$, as claimed. 
%%\equ(6.1) is an immediate consequence of this. 
We shall write $\ga$ for the common value of the $\ga(K)$.
\newline
{\bf Step 3.} In this step we deduce the almost sure convergence of
$(1/t) \log N_t(*)$. As pointed out after \equ(2), this will prove \equ(1).

We first show that $[K2^j]^{-1}Z_{K2^j}$ converges almost surely as $j \to
\infty$ for any fixed positive integer $K$. The limit turns out to be
independent of $K$.
Recall that $N_t(z(s))$ is independent of all $X(v)$ with $\|v\| \le s$
and has the same distribution as $N_t$. 
\comment
Also
$$
|\Psi(s,t)| \le t(\log d) I[N_s = 0] + s(\log d) I[N_{s,t} = 0]
\teq(11)
$$
(see the
lines preceding \equ(7.1). 
\endcomment
We now follow the second moment calculations of \cite{H} 
or \cite{SW}. We obtain from \equ(8)
$$
\frac {EZ^2_{K2^{j+1}}}{[K2^{j+1}]^2} \le 
\frac 12 \frac {EZ^2_K2^j}{[K2^j]^2} + \frac 12 \frac{[EZ_{K2^j}]^2}{[K2^j]^2}
+4\frac {\sqrt{EZ^2_{K2^j}}\sqrt{ E \Psi^2(s,t)}}{[K2^{j+1}]^2} + 
\frac {E \Psi^2(K2^j,K2^j)}{[K2^{j+1}]^2}.
\teq(12)
$$
It follows from \equ(7.7), \equ(7.1) and \equ(10) that 
$$
E\Psi^2(K2^j,K2^j) \le 2[K2^j \log (2d)]^2C_1 \exp[-C_2K2^j].
\teq(13)
$$
By subtracting $[K2^{j+1}]^{-2} [E Z_{K2^{j+1}}]^2$ from both sides of
\equ(12) and using the bound in \equ(13) we now obtain for a suitable
constant $C_3 < \infty$
$$
\text{Var}\Big[\frac {Z_{K2^{j+1}}}{K2^{j+1}}\Big] \le \frac 12 
\text{Var}\Big[\frac {Z_{K2^j}}{K2^j}\Big] +
\frac{[EZ_{K2^j}]^2}{[K2^j]^2} -
\frac {[EZ_{K2^{j+1}}]^2}{[K2^{j+1}]^2} + C_3\exp[-C_2K2^{j-1}].
\teq(14)
$$
Finally, summation of \equ(14) from $j=0$ to $j=J$  and simple algebraic
manipulations yield
$$
\frac 12 \sum_{j=0}^J\text{Var}\Big[\frac {Z_{K2^j}}{K2^j}\Big]  \le
\text{Var} \Big[\frac {Z_K}{K}\Big] + \frac {[EZ_K]^2}{K^2} + C_3\sum_{j=0}^J
C_3\exp[-C_2K2^{j-1}].
$$
Since this holds for any $J < \infty$, it follows  
$$
\sum_{j=0}^\infty\text{Var}\Big[\frac {Z_{K2^j}}{K2^j}\Big] < \infty,
$$
and then by Chebychev's inequality and Borel-Cantelli 
$$
\frac {Z_{K2^j} - EZ_{K2^j}}{K2^j} \to 0\; (j \to \infty) \text{ a.s.}
$$
Combined with \equ(9) and the independence of $\ga$ of K, this gives 
$$
\frac {Z_{K2^j}}{K2^j} \to \ga \;(j \to \infty) \text{ a.s.}
\teq(13.1)
$$

It remains to improve the convergence in \equ(13.1) to convergence along
all positive integers. To this end we fix a $0 < \ep < 1$
and note that \equ(13.1) implies 
$$
\frac {Z(\lfloor(1+\ep)^r\rfloor2^j)}{\lfloor (1+\ep)^r \rfloor 2^j}
\to \ga \text { for all integers $r \ge 0$ a.s.}
$$
Now, for small $\ep$ and for all large $n$ we  can find $1 \le r \le 
 \frac {2\log 2 }{\log(1+\ep)}$ and a $j$ such that
$$
\lfloor (1+\ep)^r\rfloor 2^j \le n \le \lfloor (1+\ep)^{r+1}\rfloor 2^j.
$$
For such $r$ and $j$ we can apply \equ(8) with $s+t = n, s = \lfloor
(1+\ep)^{r-1}\rfloor 2^j$ and 
$$
\frac \ep 3 \lfloor(1+\ep)^{r-1}\rfloor 2^j \le \frac \ep 2
(1+\ep)^{-3}n \le t=n-s \le 3 \ep  (1+\ep)^{r-1}2^j\le 4\ep (1+\ep)^{-1}n.
$$
By \equ(9.1) and \equ (10) we then have
outside a set of probability
$$
P\{N_s =0\} + P\{N_t = 0\} \le 2 C_1 \exp[-C_2 \frac \ep 2 (1+\ep)^{-3}n]
\teq(14.1)
$$
that
$$
\align
Z_n &\le Z(\lfloor (1+\ep)^{r-1} \rfloor 2^j) + Z(s,t) \\
&\le 
Z(\lfloor (1+\ep)^{r-1} \rfloor 2^j)+ t \log (2d) \text{ (see \equ(7.2))}\\
&\le Z(\lfloor (1+\ep)^{r-1} \rfloor 2^j) + 4 \log (2d)\ep (1+\ep)^{-1} n,
\teq(15)
\endalign
$$
and consequently also
$$
\frac {Z_n}n \le \frac {Z(\lfloor (1+\ep)^{r-1} \rfloor 2^j)}
{\lfloor (1+\ep)^r \rfloor 2^j} + 4\log (2d) \ep(1+\ep)^{-1}.
\teq(16)
$$
Since the sum over $n$ of the probabilities in \equ(14.1) converges,
\equ(16) will be almost surely valid for all large $n$. By taking
first the limsup as $n \to \infty$ and then as $\ep \downarrow 0$ we
find that
$$
\limsup_{n \to \infty} \frac {Z_n} n \le \lim_{\ep \downarrow 0}
\lim_{j \to \infty} \sup_{1 \le r\le 2\log 2/\log(1+\ep)}
\frac {Z(\lfloor (1+\ep)^{r-1} \rfloor 2^j)}
{\lfloor (1+\ep)^r \rfloor 2^j} = \ga \text{ a.s.}
$$
In almost the same way one can show that outside a set of 
negligable probabiliy
$$
\frac{Z (\lfloor (1+\ep)^{r+1} \rfloor 2^j)}{(1+\ep)^{r+1}2^j} 
\le (1+\ep)^2\frac {Z_n}n
$$
and obtain $\liminf_{n \to \infty} Z_n/n \ge \ga$.

We therefore proved that $\limn Z_n/n = \ga$ almost surely, and \equ(1)
with $\la = 2d e^{-\ga}$ is then immediate from the definition of $Z$.
\hfill $\blacksquare$
\enddemo

The main problem in these notes is to find information about $N_n(\al)$
as a function of $\al$. In particular, we want to compare
$\la(\al)$ to $\phi= \phi(\al) : = \limn [EN_n(\al)]^{1/n}$. 
Note that $\phi$ is easy to
evaluate. Indeed, there are $(2d)^n$ oriented paths of length $n$. A
given path $\pi$ of length $n$ contributes to $N_n$ if and only if
$W(\pi) \ge \al n.$ But, for any given $\pi$ of length $n$, $W(\pi)$
has the a binomial distribution with $n$ trials and success probability
$p$. Therefore
$$
EN_n(\al) = (2d)^n \sum_{k \ge \al n} \binom nk p^k (1-p)^{n-k},
\teq(18a)
$$
and if $\al \ge p$, then
$$
\phi = 2d \Big(\frac p \al \Big)^\al\Big(
\frac{1-p}{1-\al}\Big)^{1-\al}.
\teq(18)
$$
In the next section we shall prove a few facts concerning $\la$ and $\phi$; 
see also Figure 1.

\epsfverbosetrue
\epsfxsize=250pt
\epsfysize=180pt
\centerline{\epsfbox{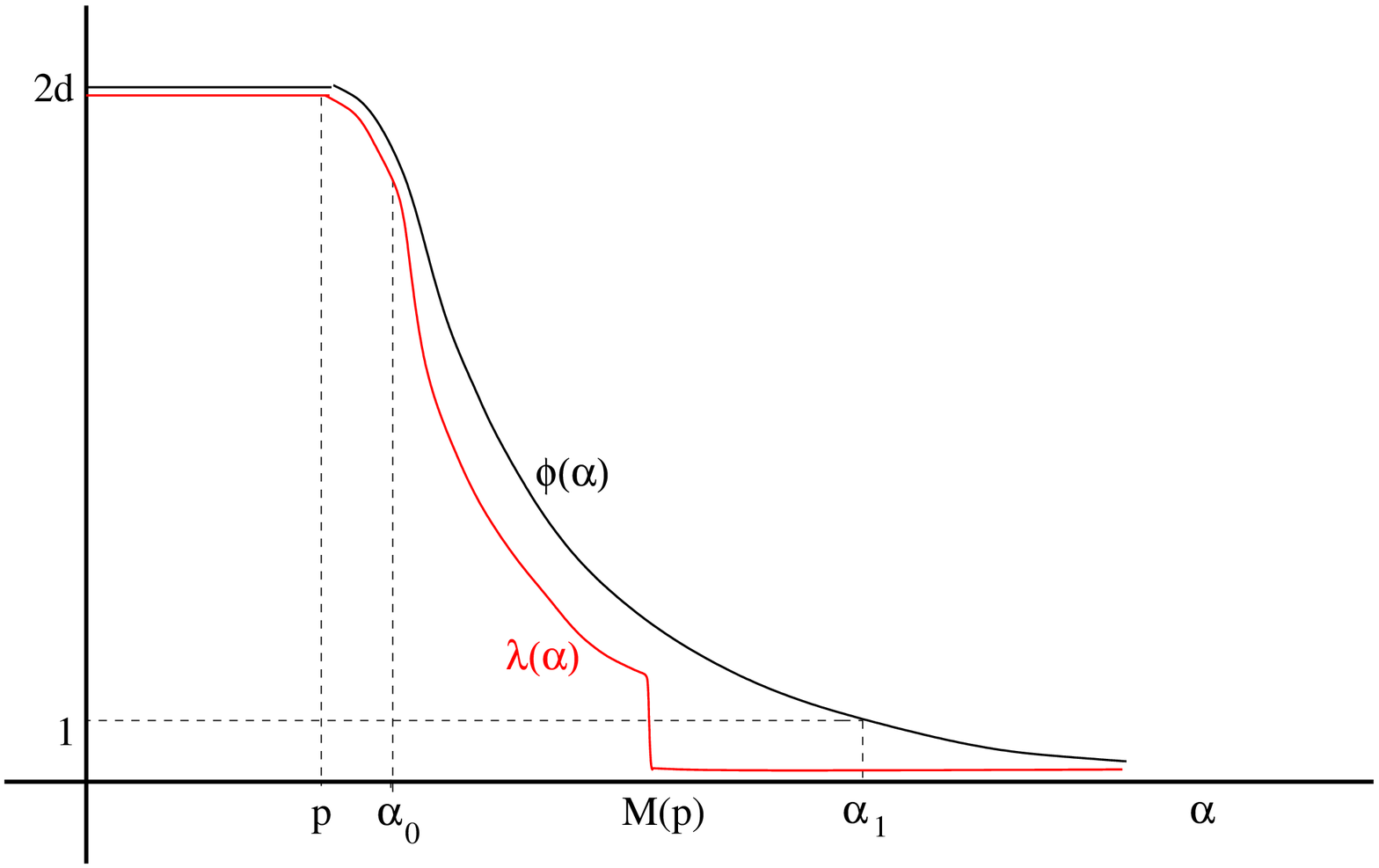}}
%%\botcaption
\noindent
Figure 1. Illustration of the graph of $\phi(\al)$ (the black curve) and of 
$\la(\al)$ (the red curve). The figure is not drawn to scale. The points 
$\al_0, \al_1$ and $M(p)$ are explained in Proposition 4, display (3.1) and 
Lemma 1, respectively.
%%\endbotcaption

\subhead
2. Properties of $\la$
\endsubhead
\numsec=2
\numfor=1
Let us first take care of the trivial region when $\al \le p$. Then
$P\{W(\pi) \ge \al n\}$ is of order 1 as $n \to \infty$ and $\phi(\al)
= 2d$. So we expect that also $\la(\al) = 2d$. The
following lemma confirms this if $\al < p$ or if $d \ge 4$ and $\al = p$.
\proclaim{Lemma 3} For $\al < p$,  or $d \ge 4$ and $\al=p$
$$
\la(\al) = \phi(\al) = 2d.
\teq(2.2)
$$
\endproclaim
\demo{Proof}
The case $d \ge 4,\al =p$, will be included in Proposition 4. 
We therefore assume throughout this proof that $\al < p$.
It is evident from the strong law of large numbers that $M(p) \ge p$, since 
$$
\limn \frac 1n W(\pi^{(n)}) = p \text{ a.s.}
$$
if $\pi^{(n)}$ is the path which moves along the first coordinate axis, i.e., 
with $\pi^{(n)}_j = (je_1, j), 0 \le j \le n$.
We therefore may assume for the rest of this proof that $\al < M(p)$.

$\phi(\al) =  2d $ for $\al \le p$ is immediate from \equ(18a) and 
the weak law of large numbers, so we
concentrate on proving $\la = 2d$. Let 
$$
R_n = \text{ (number of paths $\pi$ of length $n$ starting at 
$\bold 0$ and with } W(\pi) < \al n).
$$
Then $E R_n = (2d)^n P\{W(\pi) < \al n\}$ for any $\pi$ of length $n$ and
starting at $\bold 0$. Since $W(\pi)$ has a binomial distribution with
parameters $n,p$, and $p > \al$, 
Bernstein's inequality (\cite{CT}, Exercise 4.3.14)
shows that $P\{W(\pi) \le \al n\} \le 
 C_1 \exp[-C_2 n]$ for some constants
$C_1,C_2$ (depending on $p$ and $\al$, but not on $n$).
Consequently, $ER_n$ is exponentially small with respect to $(2d)^n$. Thus
by Markov's inequality 
$$
P\{R_n \ge \frac 12 (2d)^n\} \le \frac {ER_n}{\frac 12 (2d)^n}
$$
is also exponentially small. Hence by Borel-Cantelli, almost surely
$N_n = (2d)^n-R_n\ge \frac 12 (2d)^n$ eventually.
This, together with Lemma 2  proves 
$\la(\al) = 2d$.
\hfill $\blacksquare$
\enddemo

\comment
This further shows that if $\phi(\al) < 1$, then $N_n(\al) = 0$ for all
large $n$ almost surely, so that $\la(\al) = 0$. For small $p$ 
this will be the case for $\al \ge$ some $\al_0$ with
$$
\al_0 \sim \frac{\log d}{\log(1/p)}.
$$
We can do better, though.??????????????????????
\endcomment
The following Proposition shows that the equality $\la(\al)=\phi(\al)$
extends to $\al$ some distance beyond $p$. This is much more difficult
to prove than the preceding lemma.
\proclaim{Proposition 4} If $d \ge 4$, then 
there exists some constant $\al_0= \al_0(p) >
p$ such that 
$$
\la(\al) = \phi(\al)
\teq(2.3)
$$ 
for $\al < \al_0$. In particular $M(p) \ge \al_0$ and the limit 
$\la(\al)$ in \equ(1) exists for all $\al < \al_0$.
\endproclaim
\demo{Proof} By the proof of Lemma 3 we only have to prove \equ(2.3) 
for $p \le \al\le \al_0$ for some $\al_0 > p$. 
For the remainder of this proof a path is tacitly assumed
to have length $n$ and to start at $\bold 0$. Let 
$$
I[\pi] = \cases 1 &\text{ if }W(\pi) \ge \al n\\0
 &\text{ if }W(\pi) < \al n. \endcases
$$
Then 
$$
N_n = \sum_{\pi_0=\bold 0, |\pi| = n} I[\pi]
$$ 
We shall prove that that for suitable $\al_0 > p$
$$
EN^2_n \le C_3[EN_n]^2 \text{ for } p \le \al \le \al_0
\teq(2.5)
$$ 
for a suitable
constant $C_3 < \infty$ (independent of $n$). By Schwarz' inequality \cite
{D} this will imply 
$$
P\{N_n \ge EN_n/2\} \ge \frac 1{4 C_3}.
\teq(2.5z)
$$
In particular this will imply 
$$
M(p) \ge \limsup_{n \to \infty} \frac 1nN_n(\al) 
\ge \limsup_{n \to \infty} \frac 1{2n} EN_n(\al) = \frac 12 
\phi(\al)> 0 \text{ (see \equ(18))}.
$$
for $\al \le \al_0$. But $\limsup_{n \to \infty} \frac 1nN_n(\al') =0$ 
for $\al' > M(p)$, by definition of $M(p)$, so that $M(p) \ge \al_0$.
Lemma 2 then shows that $\la(\al) = 
\limn [N_n(\al)]^{1/n}$ exists almost surely for all $\al < \al_0$. Finally, 
\equ(2.5z) will then show that the almost sure limit of $[N_n(\al)]^{1/n}$
satisfies
$$
\la(\al) \ge \limn [EN_n(\al)]^{1/n}= \phi(\al) \text{ for }
p \le \al < \al_0.
$$

In the other direction, Markov's inequality immediately implies that always
$$
\la \le \phi.
\teq(2.1)
$$
Together these inequalities will prove \equ(2.3) and the last 
statement in the Proposition.

We turn now to the proof of \equ(2.5). Obviously
$$
\align
EN_n^2 &= \sum_{\pi'} \sum_{\pi''} E\{I[\pi']I[\pi'']\} \\ &=
\sum_{k=1}^n  \sum_{\pi'} \sum \Sb \pi'' \text{ with}\\ |\pi' \cap \pi''| = k
\endSb E\{I[\pi']I[\pi'']\}.
\teq(2.6)
\endalign
$$
Let $\{S_n'\}_{n \ge 0}$ and $\{S''_n\}_{n\ge 0}$ be two independent
simple random walks on $\cL$, both starting at $\bold 0$, and let
$T_n$ be a random variable with a binomial distribution with parameters
$n$ and $p$. Further let
$$
\rho = P\{S'_n = S''_n \text{ for some } n \ge 1\}.
$$
Then the number of pairs of paths $\pi',\pi''$ which meet at least $k$
times (not including at time 0, when both paths are at $\bold 0$) is at
most $(2d)^{2n} \rho^k$, provided $k \le n$; there are no pairs of paths
of length $n$ which meet more than $n$ times.
Let $J$ be the collection of vertices which $\pi'$ and $\pi''$ have in common
(again excluding $\bold 0$).
\comment
and write 
$$
\wt W \text{ for } \sum \Sb v \in \pi''\\ \text{ but }v \notin J
\endSb I[X(v)= e^b].
$$
\endcomment
Then, if $J$ contains exactly $k$ vertices, 
$$
\align
&P\{W(\pi'') \ge \al n\big| X(v) \text { for }v \in \pi''\}\\
& = P\Big\{\sum \Sb v \in \pi''\\ \text{ but }v \notin J \endSb  
I[X(v)= e^b] \ge \al n - \sum_{v \in J} I[X(v) = e^b]
\big| X(v) \text{ for } v \in J\Big\}\\
& \le P\Big\{\sum \Sb v \in \pi''\\ \text{ but }v \notin J \endSb  
I[X(v)= e^b] \ge \al n -k\Big\}\\
&= P\{T_{n-k} \ge \al n-k\}.
\endalign
$$
Consequently, if $|\pi' \cap \pi''| = k$, then
$$
\align
E\{I[\pi']I[\pi'']\}&= P\{W(\pi') \ge \al n\}P\{W(\pi'') \ge \al n\big|
W(\pi') \ge \al n\}\\
& \le P\{T_n \ge \al n\} P\{T_{n-k} \ge \al n - k\}\\
& \le P\{T_n \ge \al n\}.
\endalign
$$

We substitute these bounds in \equ(2.6). We then see that the right
hand side of \equ(2.6) is for any $0 < \be \le 1$ at most
$$
\align
&\sum_{1\le k \le \be n} (2d)^{2n} \rho^k P\{T_n \ge \al n\}P\{T_{n-k}
\ge \al n -k\} + (2d)^{2n}\rho^{\be n} P\{T_n \ge \al n\}\\
& \le \big[(2d)^n P\{T_n \ge \al n\}\big]^2 
\Big[\sum_{1 \le k \le \be n} \rho^k
\frac{P\{T_{n-k} \ge \al n -k\}}{P\{T_n \ge \al n\}}
+ \rho^{\be n} \frac 1{P\{T_n \ge \al n\}}\Big]\\
& = [E N_n]^2 \Big[\sum_{1 \le k \le \be n} \rho^k
\frac{P\{T_{n-k} \ge \al n -k\}}{P\{T_n \ge \al n\}}
+ \rho^{\be n} \frac 1{P\{T_n \ge \al n\}}\Big].
\teq(2.10)
\endalign
$$

Note that $\rho$ depends on $p$ and $d$ only, so is a constant $< 1$ for our
purposes here. Moreover, by \equ(18),
for any given $\al_0 > p$, it will be the case that for all $p \le
\al \le \al_0$ 
$$
\limn [P\{T_n \ge \al n\}]^{1/n} \ge \limn [P\{T_n \ge \al_0 n\}]^{1/n} 
=(2d)^{-1} \phi(\al_0) = \Big(\frac p {\al_0} \Big)^{\al_0}
\Big(\frac{1-p}{1-\al_0}\Big)^{1-\al_0}.
$$
Therefore, for any $0 \le \be \le 1$, we can choose $\al_0= \al_0(\be) > 0$ 
so close to $p$ that $\rho^{\be n} [P\{T_n \ge \al n\}]^{-1}$ is exponentially
small, uniformly in $p \le \al \le \al_0$. In other words, the second term in
the right hand side of \equ(2.10) can be taken care of by taking $\al_0-p
> 0$ small, after we have picked $\be$. Thus, to prove \equ(2.3) it
suffices to show that we can pick $\al_0> p$ and $\be > 0$ so small that
$$
\sum_{1 \le k \le \be n} \rho^k
\frac{P\{T_{n-k} \ge \al n -k\}}{P\{T_n \ge \al n\}} \le C_3-1,
\teq(2.15)
$$
uniformly for $p \le \al \le \al_0$.
%%each $\al \in [p, \al_0]$. 
Without loss of generality we 
take $\be < p$ so that $\be < \al$.

To prove \equ(2.15) we start from 
$$
\frac{P\{T_{n-k} \ge \al n - k\}}{P\{T_n \ge \al n\}}
 \le \frac {P\{T_n \ge \al n - k\}  }{P\{T_n \ge \al n \}}
=\prod_{j=1}^k \frac{P\{T_n \ge \al n -j\}}{P\{T_n \ge \al n -j+1\}}.
\teq(2.16)
$$
In addition, if for simplicity we write $\al n -j$ for $\lceil \al n -j
\rceil$, we shall use
$$
\align
P\{T_n \ge \al n - j\} &= \sum_{r = \al n -j} ^n\binom n r p^r(1-p)^{n-r}\\
&= n \binom {n-1} {\al n -j-1}\int_0^p x^{\al n-j-1} (1-x)^{(1-\al)n+j} dx,
\endalign
$$
and
$$
\frac{P\{T_n \ge \al n -j\}}{P\{T_n \ge \al n -j+1\}}
= \frac {\al n -j}{(1-\al)n+j} \frac 
{\int_0^p x^{\al n-j-1} (1-x)^{(1-\al)n+j} dx}
{\int_0^p x^{\al n-j} (1-x)^{(1-\al)n+j-1} dx}.
\teq(2.11)
$$
We want to show that the ratio here is close to 1 uniformly in 
$\al \in [p, \al_0]$ when $\al_0$ is close to $p$, 
and $1 \le j \le k \le \be n$ with $\be$ small.
We first show that we may replace the integrals over the
interval $[0,p]$ here, by integrals over $[p-\ep,1]$ for any fixed
(but sufficiently small) $\ep > 0$, without influence on the 
asymptotic behavior of 
the right hand side in \equ(2.11). To be more precise set 
$$
A= \al -\frac {j+1} n, \; B= (1-\al) + \frac jn \text{ and } f(x) = f(x;j,n)
= x^A(1-x)^B,
$$ 
so that $x^{\al n-j-1} (1-x)^{(1-\al)n+j} = f^n(x;j,n)$. Now 
$$
f'(x;j,n)= \big[\frac Ax - \frac B{1-x}\big ]f(x).
$$
One sees from this that $f(x)$ is 
strictly  increasing in $[0,A/(A+B)]$ and strictly decreasing in
$[A/(A+B),1]$. In
particular, if $ 0 \le \al -p \le \ep/4$ and $j \le \be n$ with $0 \le
\be \le (\ep /8)\land p$, then $\max_x f(x)$ is achieved at the single point
$$
x_0 := \frac A{A+B} \in [(\al - \be)-1/(n-1), \al] \subset [p-\ep/2, p+\ep/4]
$$
(provided $n \ge 1+\ep/4$) and consequently
$$
\int_0^{p-\ep} f^n(x) dx \le (p-\ep) f^n(p-\ep;n,j), 
\teq(2.12)
$$
while
$$
\int_0^p f^n(x) dx \ge \int_{p-3\ep/4}^{p-\ep/2} f^n(x) dx \ge
\frac \ep 4 f^n(p-3\ep/4;n,j).
\teq(2.13)
$$
Finally,
$$
f''(x) = -\frac A{x^2} - \frac B{(1-x)^2} \le -(A+B ) < 0
$$
so that, by Rolle's theorem,
$$
f'(x) \ge f'(p-3\ep/4) \ge f'(p-\ep/2) + (A+B) \frac \ep 4 \ge C > 0 
\text{ for $x \le p-3\ep/4$},
$$
and some constant $C = C(\ep) > 0$, independent of $\al$ and $n$.
Also, again by Rolle's theorem,
$$
f(x) \le f(p-3\ep/4) - \frac \ep 4 f'(p-3\ep/4) \le f(p -3\ep/4) - C
\frac \ep 4 \text{ for } x \le p-\ep.
$$
We combine this result with \equ(2.12) and \equ(2.13) to obtain that
$$
\align
&\frac{\int_0^{p-\ep} x^{\al n-j-1} (1-x)^{(1-\al)n+j} dx}
{\int_0^p x^{\al n-j-1} (1-x)^{(1-\al)n+j} dx}
\le \frac{4p}{\ep} \Big[\frac {f(p-\ep)}{f(p-3\ep /4)}\Big]^n\\
&\le \frac {4p}\ep \Big[1- \frac {C\ep}{4f(p-3\ep/4)}\Big]^n
\to 0\text { as } n \to \infty.
\endalign
$$
In fact, since $f(x) \le 1$, this convergence is uniform in 
$\al \in [p,\al_0]$ for $\al_0$ sufficiently close to $p$ and 
$\be$ sufficintly small.
This shows that replacement of the integral over $x  \in [0,p]$
in the numerator of the right hand side of \equ(2.11) by the same
integral over $[p-\ep,1]$, only does not change the right hand side of
\equ(2.11) much for large $n$. On the other hand, the right hand side 
of(2.9) can only increase if we
replace the integral in the denominator by the integral over $[p-\ep,p]$.
It follows that for small $\ep > 0$ and $p \le \al_0 \le p+\ep/4, 0 < \be
\le (\ep/8) \land p$, the right hand side of \equ(2.11) is 
for $p \le \al \le \al_0$ and all large $n$ at most
$$
\align
&(1+\ep)\frac {\al}{1-\al + \be} 
\frac {\int_{p-\ep}^p x^{\al n-j-1} (1-x)^{(1-\al)n+j} dx}
{\int_{p-\ep}^p x^{\al n-j} (1-x)^{(1-\al)n+j-1} dx}\\
&\le (1+\ep)\frac {p+\ep/4}{1-p-\ep/4}\cdot \frac {1-p+\ep}{p-\ep}.
\endalign
$$
Here we used that the integrand in the numerator is at most a factor
$(1-x)/x \le (1-p+\ep)/(p-\ep)$ times the integrand in the denominator for $x
\in [p-\ep,p]$. There are similar lower bounds for \equ(2.11), but we
shall not pursue these because we do not need them.

The preceding estimates show that we can choose $\ep_0 > 0$ and $\al_0 > p$
such that for $\al \in [p,\al_0]$ and all large $n$ for all 
$1 \le j \le k \le \be n$, $\rho$ times the right hand side of
\equ(2.11) is less than $1-\ep_0$ (recall that $\rho < 1$). 
The inequality \equ(2.16) then shows 
$$
\rho^k\frac{P\{T_{n-k} \ge \al n - k\}}{P\{T_n \ge \al n\}} \le
[1-\ep_0]^k \text{ for } k \le \be n,
$$
and hence also proves \equ(2.15) with $C_3 = [\ep_0]^{-1}  + 1$.
\hfill $\blacksquare$
\enddemo
\proclaim{Corollary} For $d \ge 4$ and all $0<p<1$ it holds
$$
M(p) > p.
\teq(3.6)
$$
\endproclaim
\demo{Proof}This is immediate from Proposition 4 and the fact 
that $\al_0 > p$ in this Proposition.
\enddemo
\subhead
3. Behavior of $\la(\al)$ for ``large'' $\al$ 
\endsubhead
\numsec=3
\numfor=1 
The last Proposition gives the behavior of $\la$ for ``small''
$\al$, that is, from $\al = 0$ to a little beyond $p$.
In this section we shall look at the behavior of $\la(\al)$ when 
$\la(\al)$ is small, which corresponds to large $\al$.

It is well known that on the regular $(2d)$-ary tree 
(in which each vertex has degree $2d$) it holds
$$
M(p) = \sup\{ \al: \phi(\al) > 1\}
$$ 
(see \cite {B}, Formula (3.4)). One can also 
use a branching random walk proof to show that on such a rooted
regular tree, oriented away from the root, 
for $\al$ such that
$\phi(\al)> 1$, it holds $\la(\al) = \phi(\al)$. As we shall demonstrate
soon, this is not the case
for walks on $\cL =\Bbb Z^d \times \Bbb Z_+$. 

If $\al$ is such that $\phi(\al) < 1$, then, by the definition of
$\phi$, $EN_n(\al)$ tends to 0 exponentially fast, so almost surely
$N_n(\al) = 0$ eventually. Of course $\la(\al)=0$ in this case. If $p$
is small, this case applies for $\al >$ some $\al_1$ with
$$
\al_1 \sim \frac{\log (2d)}{\log(1/p)}.
\teq(3.1)
$$
We can do better, though. By definition of $M(p)$, if $\al > M(p)$,
then $N_n(\al) =0$ for large $n$. Thus 
$$
\la(\al) = 0 \text{ if } \al > M(p).
\teq(3.2)
$$
But it is shown in \cite{L} that there exist constants $C_1, C_2 \in
(0, \infty)$ such that
$$
C_1p^{1/(d+1)} \le M(p) \le C_2 p^{1/(d+1)}.
\teq(3.3)
$$
Thus, by \equ(3.2), for small $p$ it holds
$$
\la(\al) = 0 \text{ if } \al > C_2p^{1/(d+1)}.
\teq(3.4)
$$
Clearly this improves \equ(3.1) for small $p$; it shows that $\la(\al)$
is still zero for smaller values of $\al$ than indicated by \equ(3.1).
We shall next show that \equ(3.2) is best possible in the following
sense.
\proclaim{Proposition 5} For $d \ge 4$ and each $p \in (0,1)$ it holds
$$
\la(\al) > 1 \text{ for all } \al < M(p).
\teq(3.5)
$$
\endproclaim
\comment
{\bf Remark:} $p_0$ merely needs to satisfy 
$$
M(p) > p \text{ for all }0 < p \le p_0.
\teq(3.6)
$$ 
Such a $p_0> 0$ exists by \equ(3.3).

It is possibly that any $p_0 < p_c(\Bbb Z^d,\text{ oriented,
site)}$ will work. What is required of $p_0$ is merely that 
there exist constants $C_i \in (0, \infty)$ such that
$$
P_p\{\exists \text{ a path $\pi$ of length $n$ with $W(\pi) \ge
(1-C_4)n$} \} \le C_5\exp[-C_6n]
\teq(3.6)
$$
for any $p \le p_0$.
\endcomment
\demo{Proof} For $\al \le p$, \equ(2.2) already shows that 
$\la(\al) = 2d > 1$. For the remainder of this proof we therefore take 
$\al > p$. As before it is tacitly assumed that all paths in this
proof start at $\bold 0$.
\comment
 We shall also assume that \equ(3.6) holds
for $0< p \le p_0$.
Fix $0< p \le p_0, p < \al < M(p)$ and let $0 < \eta \le 1/4$. 
\endcomment
Fix $\eta \in (0, 1/4)$ and define $\wt \al = [\al+M(p) ]/2$, 
so that $p < \wt \al < M(p)$ by \equ(3.6).
By Theorem 2 in \cite{GK} there then exists an 
$M_0 < \infty$ such that with probability at least $(1-\eta)$ there 
exists for each $n \ge M_0$ a path $\wt \pi$
starting at $\bold 0$ and of length $n$ which has $W(\wt \pi) \ge \wt \al n$.
Now fix $n \ge M_0$ and let $\wt \pi = (\bold 0 = \wt \pi_0, \wt
\pi_1, \dots \wt \pi_n)$ be a path with the above properties. Assume
that for a certain $k \le n-2$
$$
e(i_{k+1}) := \wt \pi_{k+2} -\wt \pi_{k+1} \ne e(i_k) := \wt \pi_{k+1} - \wt
\pi_k.
\teq(3.7)
$$
%%and 
%%$$
%%X(\wt \pi_{k+1}) = 0
%%\teq(3.8)
%%$$
We can then interchange the two steps $e(i_k)$ and $e(i_{k+1})$ to get
the new path
$$
\wh \pi =(\bold 0, \wt \pi_1, \dots, \wt \pi_k, \wt \pi_k +
e(i_{k+1}), \wt \pi_k + e(i_{k+1}) + e(i_k)=\wt \pi_{k+2}, \wt
\pi_{k+3}, \dots, \wt \pi_n).
$$
This path differs only in its point at time $k+1$ from $\wt\pi$,
so that
$$
W(\wh \pi)-W(\wt \pi) \ge - X(\wt \pi_{k+1}) \ge -1.
\teq(3.8)
$$
However $\wh \pi$ will still be selfavoiding, since $\wt \pi$
does not visit $\wh \pi_{k+1}$,  because $\|\wh \pi_{k+1}\| = \|\wt
\pi_{k+1}\| = k+1$ and $\wt \pi$ can  visit only one point with
$\l_1$-norm $k+1$. 
\comment
Morever, the point $\wt \pi_{k+1}$
contributed only weight $0$ to $W(\wt \pi)$ (see \equ(3.8)) so that in
any case
$$
W(\wh \pi)-W(\wt \pi) = W(\wh \pi) \ge 0
$$
and hence
$$
W(\wh \pi) \ge W(\wt \pi) \ge  \al n.
$$
Thus,$\wh \pi$ also contributes one to the count $N_n(\al)$.
\endcomment
If there
are $m$ values of $k$ for which \equ(3.7) holds, then we can 
interchange two successive steps as described above or not at at least $m/2$
places such that these interchanges do not interfere with each other
(say, at any subset of the even $k$'s which satisfy \equ(3.7)). 
This yields at least $2^{m/2}$ paths with weight $W \ge \wt \al
n- m/2$. In other words, $[N_n(\wt \al - m/(2n))]^{1/n} \ge 2^{m/(2n)}$
in this case. If we take 
$$
0 < \liminf_{n \to \infty} m/(2n) \le \limsup_{n \to \infty}m/(2n) \le \wt \al
- \al
$$
then this method results 
in 
$$
\liminf_{n \to \infty} [N_n(\al)]^{1/n} \ge \exp[\liminf_{n \to
\infty} m/(2n) \log 2] > 1.
$$ 

In view of the preceding paragraph and the fact that $\lim_n
[N_n(\al)]^{1/n}$ exists, it suffices for the  
proposition that there is for all large $n$ 
at least a probability $\eta$ that there is
a path $\wt \pi$ of length $n$ and $W(\wt \pi) \ge \wt \al n$ and for
which \equ(3.7) holds for at least $C_3 n$ values of $k$
(with $C_3 > 0$ and independent of $n$ and $\wt \pi$). In this case we
may take $m = (C_3 \land (\wt \al-\al))n$ in the preceding argument.

Let us now make sure that we can find $\wt \pi$ so that $W(\wt \pi)
\ge \wt \al n$ and such that \equ(3.7) holds for many $k$. We shall bound 
the probability that no such path exists. This last probability
is, for $n \ge M_0$, bounded by
$$
\align
&P\{\text{there is no path $\pi$ of length $n$ with $W(\pi) \ge \wt \al n$}\}\\
&+P\{\text{there exists a path $\wt \pi$ of length $n$ with 
$W(\wt \pi) \ge \wt \al n$ but}\\
&\phantom{MMMMMMMMM}
\text{fewer than $C_3n$ values of $k$ for which \equ(3.7) holds}\}\\
&\le \eta + \text{(number of paths $\wt \pi$ for which \equ(3.7) holds
for no more}\\
&\phantom{MMMMMMMMM}\text{than $C_3 n$ values of $k$})
\,P\big\{\sum_{i=1}^n Y_i \ge \wt \al n\big\},
\teq(3.10)
\endalign
$$
where the $Y_i$ are i.i.d., each with the distribution $P\{Y_i =1\} =
1-P\{Y_i =0\} = p$.
But any path $\wt \pi$ of length $n$ is determined by  
the values of the $k$ for which \equ(3.7) holds as well as the values
of the corresponding $e(i_{k+1})$, and also $\wt \pi_1$. Indeed this
gives the places at which the direction of the steps of $\wt \pi$
changes and the value of this direction immediately after the change
(plus the starting direction). The number of paths for which \equ(3.7)
holds for no more than $C_3n$ values of $k$ 
and the number of choices for the directions right after the $k_i$ and
at time 0 is at most 
$$
2d\sum_{j \le C_3n} \binom n j (2d-1)^j \le C_4 \exp[nC_3 \log (1/C_3)
+n C_3\log (2d-1)]
$$ 
for small $C_3$.
But $\wt\al > p$, and by simple exponential bounds for the binomial
distribution (e.g., Bernstein's inequality in \cite{CT}, Exercise 4.3.14)
$$
P\{\sum_{i=1}^n Y_i \ge \wt \al n\} \le C_5\exp[-C_6n] 
$$
for some constants $0 < C_5(p, \wt \al), C_6(p, \wt \al)< \infty$.
Thus the right hand side of \equ(3.10) is bounded by
$$
\eta +  C_4C_5 \exp[nC_3 \log (1/C_3)
+n C_3\log (2d-1)- nC_6].
$$ 
Since $C_6$ is independent of $C_3$,
we can choose $C_3> 0$ so small that this
expression is at most $2\eta$ for large $n$. The complementary
probability is then
$$
\align
&P\{\text{there exist a path $\wt \pi$ of length $n$ with $W(\wt \pi)
\ge \wt \al n$},\\
&\phantom{MMMMM}\text{and all such paths have at least $C_3n$ values of $k$ for
which \equ(3.7) holds}\}\\
& \ge 1- 2\eta > \eta
\text{ (recall }\eta \le 1/4). 
\tag "$\blacksquare$"
\endalign
$$
\enddemo

\medskip

\noindent {\bf References.}

\medskip 

\noindent [B] Biggins, J.D. (1977) Chernoff's theorem in the 
branching random walk. J. Appl. Probab. {\bf 14}, 630--636.

\smallskip 

\noindent [CT] Chow, Y. S.; Teicher, H. (1988), Probability theory. 
Independence,
Interchangeability, Martingales. Second edition. Springer Texts in
Statistics. Springer-Verlag, 1988.

\smallskip 

\noindent [CPV] Comets, F. Popov, S, Vashkovskaya, M. (2007) Private
communication.

\smallskip

\noindent [CMS] Cranston, M., Mountford, T. S., Shiga, T. (2005) Lyapounov 
exponents for the parabolic Anderson model. Acta Math. Univ. Comenian. 
{\bf 71}, 163--188. 

\smallskip 

\noindent [GK] Gandolfi, A.; Kesten, H., (1994) Greedy lattice animals II. 
Linear growth.  Ann. Appl. Probab. {\bf 4}, 76--107.

\smallskip 

\noindent [H] Hammersley, J. M. (1974) Postulates for subadditive processes.
Ann. Probab. {\bf 2}, 652--680. 

\smallskip 

\noindent [D] Durrett, R. (1996) Probability: Theory and Examples. 
Second edition. Duxbury Press, Belmont, CA, 1996.

\smallskip 

\noindent [L] Lee, S. (1994). A note on greedy lattice animals. 
Ph.D. dissertation, Cornell Univercity, Ithaca, NY.

\smallskip 

\noindent [Li] Liggett, T. M. (1985) Interacting Particle Systems. 
Springer-Verlag, 1985.

\smallskip 

\noindent [SW] Smythe, R. T.; Wierman, J. C. (1978) First-passage 
Percolation on the Square Lattice. Lecture Notes in Mathematics, vol. 671. 
Springer-Verlag, 1978.

\enddocument